\documentclass[11pt]{article}
\usepackage{amsmath}
\usepackage{amscd}
\usepackage{amsfonts}
\usepackage{amssymb}
\usepackage{amsthm}

\newcommand{\rar}{\rightarrow}
\newcommand{\calg}{\mathcal}

\newcommand{\strc}{\mathcal{O}_{X}}

\newcommand{\dcat}[1]{\text{D}^{b}_{\text{coh}}(#1)}

\newcommand{\hkr}{\text{I}_{HKR}}

\newcommand{\cale}{\mathcal{E}}

\newtheorem{thm}{Theorem}
\newtheorem*{thmC}{Explicit Cardy condition}
\newtheorem{prop}{Proposition}

\title{A generalized Hirzebruch Riemann-Roch theorem.}

\author{Ajay C. Ramadoss}

\begin{document}

\maketitle

\begin{abstract}

This short note proves a generalization of the Hirzebruch
Riemann-Roch theorem equivalent to the Cardy condition described in
[1]. This is done using an earlier result [4] that explicitly
describes what the Mukai pairing in [1] descends to in Hodge
cohomology via the Hochschild-Kostant-Rosenberg map twisted by the
root Todd genus.

\end{abstract}

\section{Statement of the generalized Riemann-Roch theorem.}

Let $X$ be a smooth proper scheme of dimension $n$ over a field
$\mathbb K$ of characteristic $0$. Let $\mathcal{E}$ and
$\mathcal{F}$ be elements of the bounded derived category
$\text{D}^{b}_{\text{coh}}(X)$ of coherent $\strc$-modules on $X$.
Let $v$ and $w$ be elements of $\text{End}_{\dcat{X}}(\mathcal{E})$
and $\text{End}_{\dcat{X}}(\calg F)$ respectively. Let
$\text{T}_{v,w}$ denote the endomorphism $$a \mapsto w \circ a \circ
v $$ of $\text{RHom}_X(\mathcal{E},\mathcal{F})$. Let
$\text{at}(\calg E) \in \text{Hom}_{\dcat{X}}(\cale,\cale \otimes
\Omega[1])$ denote the Atiyah class of $\calg E$. Let
$\text{ch}_v(\mathcal{E})$ denote the "twisted" Chern character
$$\text{Tr}_{\cale}(\text{exp}(\text{at}(\mathcal{E})) \circ v) \in
\text{Hom}_{\text{D}^{b}_{\text{coh}}(X)}(\mathcal{O}_X,\oplus_i
\Omega^i_X[i]) \simeq \oplus_i \text{H}^i(X,\Omega^i_X)\text{ . }$$
Note that if $v = \text{id}_{\mathcal E}$ then $\text{ch}_v(\mathcal
E)= \text{ch}(\mathcal E)$, the Chern character of $\mathcal E$. Let
$K$ be the involution on $\oplus_{p,q} \text{H}^{q}(X,\Omega^{p})$
that acts on the summand $\text{H}^{q}(X,\Omega^{p})$ by
multiplication with ${(-1)}^q$. If $x$ in an element of
$\oplus_{p,q} \text{H}^{q}(X,\Omega^{p})$, we shall denote $K(x)$ by
$x^*$. For any endomorphism $\text{T}$ of
$\text{RHom}_X(\mathcal{E},\mathcal{F})$, $\text{str}(\text{T})$
shall denote the alternated trace of $\text{T}$. If $f:X \rar Y$ is
a morphism of smooth proper schemes, $f_*$,$f^*$ etc shall denote
the corresponding derived functors unless explicitly mentioned
otherwise. $\int_X$ shall denote the linear functional on
$\oplus_{p,q} \text{H}^{q}(X,\Omega^{p})$ that coincides with the
Serre duality trace on $\text{H}^n(X,\Omega^n_X)$ and vanishes on
all other direct summands. We have the following generalization of
the Hirzebruch Riemann-Roch theorem.
\\

\begin{thmC}
$$\text{str}(\text{T}_{v,w}) = \int_X \text{ch}_v(\mathcal{E})^*
\text{ch}_w(\mathcal{F}) \text{td}(X) \text{ . }$$
\end{thmC}

\textbf{Remark 1:} Note that when $v = \text{id}_{\mathcal E}$ and
when $w=\text{id}_{\mathcal F}$ the above statement amounts to the
Hirzebruch Riemann-Roch theorem

$$\chi({\mathcal E},{\mathcal F}) = \int_X \text{ch}(\mathcal
E)^*\text{ch}(\mathcal F)\text{td}(X) \text{ . }$$

 \textbf{Remark 2:} The proof of the above theorem will also enable
us to see that the above theorem is in fact equivalent to the Cardy
condition in [1] (see Theorem 7.9 of [1]). The more classical
(compared to Theorem 7.9 of [1]) statement of the Cardy condition
given here would be useful for those interested in an explicit
version of the Cardy condition for computational purposes. Such
computations are related to understanding string propagation between
D-Branes with twisted boundary conditions in the situation where $X$
is Calabi-Yau and the category of D-Branes is equivalent to
$\text{D}^{b}_{\text{coh}}(X)$.
An analogous computation in a different set-up seems to have been discussed at least implicitly in [5]. \\

\subsection*{Acknowledgements}

I am deeply grateful to Prof. Boris Tsygan for informing me about
this problem and providing me with a crucial insight. It was Prof.
Tsygan who told me that conjecturally, there should be a
Riemann-Roch theorem amounting to the Cardy condition relating
$\text{str}(\text{T}_{v,w})$ to $\text{ch}_v(\mathcal{E})$,
$\text{ch}_w(\mathcal{F})$ and the Todd genus. Prof. Tsygan also
informed me that he was asked this question by Prof. David Kazhdan
and Prof. Anton Kapustin. My heartfelt thanks to Prof. Kazhdan and
Prof. Kapustin as well.

\section{The Mukai pairing and the Cardy condition.}

Let $\Delta:X \rar X \times X$ be the diagonal embedding. Let
$\Delta_!$ denote the left adjoint of $\Delta^*$. Recall that the
$i$-th Hochschild homology $\text{HH}_i(X)$ is the space
$\text{Hom}^i_{\text{D}^{b}_{\text{coh}}(X \times X)}(\Delta_!
\strc, \Delta_* \strc)$. In his paper [1], Caldararu defined a
(nondegenerate) Mukai pairing
$$\langle \text{  },\text{  }\rangle_M :\text{HH}_i(X) \otimes
HH_{-i}(X) \rar \mathbb K \text{  } \forall \text{  } i \text{ . }$$

In [1], we have a map
$$\iota^{\mathcal{H}}:
\text{Hom}_{\dcat{X}}(\mathcal{H},\mathcal{H}) \rar \text{HH}_{0}(X)
 = \text{Hom}_{\dcat{X \times X}}(\Delta_! \strc,\Delta_* \strc)
$$ for any $\mathcal{H} \in \dcat{X}$. Also, $\Delta_! \strc \simeq
\Delta_* S_X^{-1}$ where $S_X$ is the shifted line bundle tensoring
with which yields the Serre duality functor on $\dcat{X}$. Let $\nu
\in \text{Hom}_{\dcat{X \times X}}(\Delta_* \strc,\Delta_* S_X)$.
For any smooth proper scheme $Y$ and any $\calg V \in \dcat{Y}$, let
$\text{Tr}_Y:\text{Hom}_{\dcat{Y}}(\calg V,S_Y \otimes \calg V) \rar
\mathbb K$ denote the Serre duality trace. If $u \in
\text{Hom}_{\dcat{X}}(\calg H,\calg H)$, then

\begin{prop}(see the definition of $\iota^{\mathcal H}$ in Section 6.3 of [1])
$$ \text{Tr}_{X \times X}(\nu \circ \iota^{\calg H}(u)) =
\text{Tr}_X (\pi_{2*}(\pi_1^* \calg H \otimes \nu) \circ u) \text{ .
}$$

\end{prop}

The following statement of the Cardy condition appears as Theorem
7.9 in [1]. Let $\calg E$,$\calg F$,$v$,$w$ be as in the previous
subsection.

\begin{thm}
$$\text{str}(T_{v,w})= \langle \iota^{\cale}(v),\iota^{\calg
F}(w)\rangle_M \text{ . }$$
\end{thm}

Recall that the Adjunction $\Delta_! \dashv \Delta^*$ identifies
$\text{Hom}^i_{\dcat{X \times X}}(\Delta_! \strc, \Delta_* \strc)$
with $\text{Hom}^i_{\dcat{X}}(\strc,\Delta^*\Delta_* \strc)$.
Therefore, the Hochschild-Kostant-Rosenberg map
$$\hkr:\Delta^*\Delta_* \strc \rar \oplus_j \Omega^j[j] $$ induces
maps
$$\hkr:HH_i(X) \rar \oplus_j \text{H}^{j-i}(X,\Omega^j_X) \text{ .
}$$ The following theorem appears implicitly in Markarian's work [3]
and in the form stated below in an earlier paper of this author [4].
Let $x \in \text{HH}_i(X)$ and let $y \in \text{HH}_{-i}(X)$. Then,

\begin{thm}
$$\langle x,y\rangle_M = \int_X \hkr(x)^*\hkr(y)\text{td}(X) \text{
. }$$
\end{thm}

In order to prove the explicit Cardy condition, we therefore need to
prove the following proposition.

\begin{prop}
$$\hkr(\iota^{\cale}(v)) = \text{ch}_v(\cale) \text{ . }$$
\end{prop}

Note that setting $x= \iota^{\cale}(v)$ and $y =\iota^{\calg F}(w)$
in Theorem 2 ,applying Theorem 1 and using the fact that
$\hkr(x)=\text{ch}_v(\cale)$ and that $\hkr(y)= \text{ch}_w(\calg
F)$ by Proposition 2, we obtain the explicit Cardy condition. We now prove Proposition 2 below.\\

\begin{proof} This proof is a straightforward modification of the proof of Theorem 4.5 in [2]. Remember that if
$\mathcal A, \mathcal H , \mathcal J \in \dcat{X}$, there is a trace
map (see Section 2.4 of [1])
$$\text{Tr}_{\mathcal A}: \text{Hom}_{\dcat{X}}(\mathcal A \otimes \mathcal H, \mathcal A \otimes \mathcal J) \rar \text{Hom}_{\dcat{X}}
(\mathcal H,\mathcal J) \text{ . }$$ For any $\theta \in
\text{Hom}_{\dcat{X}}(\mathcal H,\mathcal J)$, we will abuse
notation and denote the map $\text{id}_{\mathcal A} \otimes \theta
\in \text{Hom}_{\dcat{X}}(\mathcal A \otimes \mathcal H,\mathcal A \otimes \mathcal J)$ by $\mathcal A \otimes \theta$.  \\

Recall that the Adjunction $\Delta_! \dashv \Delta^*$ identifies
$\text{Hom}_{\dcat{X \times X}}(\Delta_! \strc, \Delta_*\strc)$ with
$\text{Hom}_{\dcat{X}}(\strc,\Delta^*\Delta_* \strc)$. Let
$\hat{\text{ch}}_v(\cale)$ denote the image of $\iota^{\cale}(v)$
under this identification. The desired proposition states that the
image of $\hat{\text{ch}}_v(\cale)$ under $\hkr$ is
$\text{ch}_v(\cale)$. Let $\nu$ be an arbitrary element of \\
$\text{Hom}_{\dcat{X}}(\Delta^* \Delta_* \strc, S_X)$. This
corresponds to the element
$$\bar{\nu} = \Delta_* \nu \circ \eta \in \text{Hom}_{\dcat{X \times
X}}(\Delta_* \strc,\Delta_* S_X) $$ where $\eta$ is the unit of the
adjunction $\Delta^* \dashv \Delta_*$ applied to $\Delta_* \strc$.\\

Note that $$\text{Tr}_X(\nu \circ \hat{\text{ch}}_v(\cale)) =
\text{Tr}_{X \times X}(\bar{\nu} \circ \iota^{\cale}(v)) \text{( see
the proof of Proposition 1 of [4] for instance)}$$
$$ = \text{Tr}_X(\pi_{2*}(\pi_1^*\cale \otimes \bar{\nu}) \circ v)
\text{ ( by Proposition 1 )} $$ $$ =
\text{Tr}_X(\pi_{2*}(\pi_1^*\cale \otimes (\Delta_* \nu \circ \eta))
\circ v) = \text{Tr}_X(\pi_{2*}((\pi_1^*\cale \otimes \Delta_* \nu)
\circ (\pi_1^* \cale \otimes \eta)) \circ v) $$
$$=\text{Tr}_X((\cale \otimes \nu) \circ (\pi_{2*}(\pi_1^* \cale
\otimes \eta)) \circ v) $$
$$ = \text{Tr}_X(\nu \circ \text{Tr}_{\cale}(\pi_{2*}(\pi_1^* \cale
\otimes \eta) \circ v)) \text{( by Lemma 2.4 of [1]). }$$

By the non-degeneracy of the Serre duality pairing, it follows that
$$ \hat{\text{ch}}_v(\cale) = \text{Tr}_{\cale}(\pi_{2*}(\pi_1^* \cale
\otimes \eta) \circ v) \text{ . }$$ Hence,
$$\hkr(\hat{\text{ch}}_v(\cale)) =
\hkr(\text{Tr}_{\cale}(\pi_{2*}(\pi_1^* \cale \otimes \eta) \circ v)
$$ $$ = \text{Tr}_{\cale}(\pi_{2*}(\pi_1^* \cale \otimes \Delta_*
\hkr \circ \eta) \circ v) \text{ . }$$ By Proposition 4.4 of [2],
$$\Delta_* \hkr \circ \eta = \text{exp}(\alpha) $$ where
$\alpha:\Delta_* \strc \rar \Delta_* \Omega_X[1]$ is the universal
Atiyah class. It follows that
$$\text{Tr}_{\cale}(\pi_{2*}(\pi_1^* \cale \otimes \Delta_*
\hkr \circ \eta) \circ v) = \text{ch}_v(\cale) \text{ . }$$

\end{proof}

\section*{References.}

[1] Caldararu, A., The Mukai pairing I : the Hochschild structure.
Arxiv preprint
math.AG/0308079.\\

[2]. Caldararu, A., The Mukai pairing II : the
Hochschild-Kostant-Rosenberg isomorphism. {\it Advances in
Mathematics.} \textbf{194}(2005),no.1,34-66.\\

[3] Markarian, N., Poincare-Brikhoff-Witt isomorphism, Hochschild
homology ad the Riemann-Roch theorem. MPI preprint MPI 2001-52.\\

[4] Ramadoss, A., The relative Riemann-Roch theorem from Hochschild
homology. Arxiv preprint math.AG/0603127.\\

[5] Ishikawa, H., Tani, T., Twisted boundary states and
representation of generalized fusion algebra. Arxiv preprint
arxiv:hep-th/0510242.\\

\textbf{Address:}\\
 Department of Mathematics\\
 Cornell University\\
 580 Malott Hall\\
 Ithaca, NY-14853.\\

 \textbf{Email:}ajaycr@math.cornell.edu

\end{document}